\documentclass{article}
\usepackage[margin=35mm]{geometry}
\usepackage{amsthm,amsmath,amssymb}
\usepackage{mathtools}

\newtheorem{theorem}{Theorem}[section]
\newtheorem{proposition}[theorem]{Proposition}

\newtheorem{lemma}[theorem]{Lemma}
\newtheorem*{theorem*}{Theorem}

\theoremstyle{remark}
\newtheorem*{remark}{Remark}
\newtheorem*{remarks}{Remarks}
\newtheorem*{notation}{Notation}

\newcommand{\tZ}{\mathtt 2}                            
\newcommand{\tL}{\mathtt 1}                            
\newcommand{\tO}{\mathtt 0}                            

\DeclareMathOperator{\e}{\mathrm{e}}                   
\DeclareMathOperator{\LandauO}{\mathcal O}             
\def\bl{\hspace{0.5pt}\underline{\hphantom{\hspace{0.6em}}}\hspace{0.5pt}}	

\usepackage{tikz}
\usepackage{pgfplots}
\usetikzlibrary{matrix,shapes,arrows,positioning,shadows,backgrounds}

\title{Collisions of digit sums in bases $2$ and $3$}
\author{
Lukas Spiegelhofer\footnote{Supported by the FWF-ANR project ArithRand (grant numbers I4945-N and ANR-20-CE91-0006).}
\\Montanuniversit\"at Leoben, Austria
}
\date{\small\textbf{Dedicated to Jean-Marc Deshouillers on the occasion of his 75th birthday}}
\begin{document}
\maketitle
\begin{abstract}
We prove a folklore conjecture concerning the sum-of-digits functions in bases two and three: there are infinitely many positive integers $n$ such that the sum of the binary digits of $n$ equals the sum of the ternary digits of $n$.
\end{abstract}
\maketitle
\renewcommand{\thefootnote}{\fnsymbol{footnote}} 
\renewcommand{\thefootnote}{\arabic{footnote}}
\section{Introduction and main result}\label{sec_introduction}
Representations of the same number $x$ in two or more multiplicatively independent integer bases apparently look very different.
This topic is far from being understood, and the relation of the base-$q_1$ and the base-$q_2$ expansion to each other is a source of difficult problems.

The base-$q$ expansion is intimately connected to powers of $q$.
In order to understand the relation of different bases $q_1$ and $q_2$ to each other better we consider, as a start, the arrangement of powers of $2$ and $3$.
Assume that the set containing all powers of two and three (with nonnegative exponents) is sorted in ascending order:
\[(a_n)_{n\geq 0}=
(1,2,3,4,8,9,16,27,32,64,81,128,243,256,512,729,1024,\ldots)
\]
(this is sequence~\texttt{A006899} in Sloane's OEIS~\cite{OEIS}).
In what manner are the powers of two and three interleaved?
Taking logarithms, we see that the answer to this question is encoded in the Sturmian word
\[w=\bigl(\lfloor (n+1)\alpha\rfloor-\lfloor n\alpha\rfloor\bigr)_{n\geq 0},\]
where $\alpha=\log 3/\log 2=\log_2(3)$, as follows:
start with $3^0=1$, append the first $w_0=1$ powers of two --- that is, the integer $2$ --- append $3^1$, then $w_1=2$ powers of two, followed by $3^2$ and $w_2=1$ powers of two, and so on.
Our question is therefore equivalent to understanding the
\emph{continued fraction expansion} of $\alpha$ (consult, for example, Berth\'e~\cite{Berthe2001} for an explanation of this connection).
However, it is not even known whether the sequence of partial quotients of $\alpha$ is bounded, that is, whether $\alpha$ is \emph{badly approximable};
any system in this sequence has yet to be found. 
The number $\alpha$ is transcendental by the Gelfond--Schneider theorem~\cite{Gelfond1934b,Gelfond1934a};
by Baker's theorem~\cite{Baker1966,Baker1967a,Baker1967b} we obtain
\[\left\lvert\frac{\log 3}{\log 2}-\frac pq\right\rvert\geq\frac{c}{q^\rho}\]
for all integers $q>0$ and $p$ and some effective positive constants $c$ and $\rho$.
More precisely, a bound for the \emph{irrationality measure} $\mu(\alpha)$ of $\alpha$, which is the infimum of $\rho$ for which there exists $c$ such that this estimate holds for all $p,q$, was given by Rhin~\cite[Equation~(8)]{Rhin1987}: we have $\mu(\alpha)\le8.616$.
Also, Wu and Wang~\cite{WuWang2014} obtained the bound $\mu(\log 3)\le5.1163051$.
Note that badly approximable numbers have irrationality measure $2$.
We would also like to mention the interesting blog entry by Tao\footnote{
\texttt{https://terrytao.wordpress.com/2011/08/21/hilberts-seventh-problem-and-powers-of-2-and-3/}} on the topic.

In view of the above problem we have to expect major difficulties when we try to mix different bases.
In this context, the following unsolved conjecture of Furstenberg~\cite{Furstenberg1970} is of interest, concerning \emph{multiplicatively independent} integer bases $p,q\geq 2$
(that is, such that $p^k\neq q^\ell$ for all $k,\ell\geq 1$):
define
\[O_a(x)\coloneqq \bigl\{a^kx \bmod 1:k\in\mathbb N\bigr\}\]
and let $\dim_H(A)$ be the Hausdorff dimension of a set $A\subseteq[0,1]$.
Then
\begin{equation}\label{eqn_furstenberg}
\dim_H\bigl(\overline{O_p(x)}\bigr)+\dim_H\bigl(\overline{O_q(x)}\bigr)\geq 1
\end{equation}
for all irrational $x\in[0,1]$.
Furstenberg's conjecture underlines the idea stated before: different bases should produce very different representations of the same number.
We note the papers~\cite{Shmerkin2019,Wu2019} for recent progress on this conjecture, and the recent preprint~\cite{AdamczewskiFaverjon2020} by Adamczewski and Faverjon, where related independence results can be found.

The related topic of studying the base-$p$ expansion of powers of $q$ is very difficult and has attracted the attention of many researchers; we note the recent preprint~\cite{KerrMeraiShparlinski2021} by Kerr, M\'erai, and Shparlinski and the references contained therein.
Erd\H{o}s~\cite{Erdos1979} conjectured that the only powers of two having no digit $\tZ$ in its ternary expansion are $1,4$, and $256$ (see also Lagarias~\cite{Lagarias2009}).
This conjecture is open, and Erd\H{o}s wrote ``[\ldots] as far as I can see, there is no method at our disposal to attack this conjecture''~\cite{Erdos1979}.
Meanwhile, there is a close connection to to Erd\H{o}s' \emph{squarefree conjecture}~\cite{Erdos1980}, stating that the central binomial coefficient $\binom{2n}n$ is never squarefree for $n\geq 5$.
The latter conjecture was proved for all large $n$ by Sárk\"ozy~\cite{Sarkozy1985}, and solved completely by Granville and Ramar\'e~\cite{GranvilleRamare1996}.
The connection between these two conjectures can be understood by considering the identities
\[\nu_2\left(\binom{2n}n\right)=s_2(n)
\quad\text{and}\quad
\nu_3\left(\binom{2n}n\right)=s_3(n)-\frac{s_3(2n)}2,  \]
where $s_q$ is the sum-of-digits function in base $q$, and $\nu_p$ is the $p$-adic valuation of an integer $\geq 1$ (with $p$ prime).
That is, $\binom{2n}n$ is divisible by the square $4$ if $n\geq 1$ is not a power of two,
and so a stronger form of the (already proved) squarefree conjecture would follow from a proof of the conjecture that
\begin{equation}\label{eqn_div_conj}
s_3(2^k)-s_3(2^{k+1})/2\geq 2\quad\text{for }k\geq 9.
\end{equation}
In fact,~\eqref{eqn_div_conj} implies $4\mid\binom{2n}n$ or $9\mid\binom{2n}n$ for each $n\geq 257$, while $\binom{512}{256}$ is divisible by neither $4$ nor $3$.
Equation~\eqref{eqn_div_conj} in turn would follow if we could prove that the integer $2^k$ contains at least two digits equal to $\tZ$ in ternary for $k\geq 9$: in this case at least two carries appear in the addition $2^k+2^k$ in ternary.
We also would like to note the recent preprint~\cite{DimitrovHowe2021} by Dimitrov and Howe on this topic.

The main objects in this paper are the sum-of-digits functions $s_2$ and $s_3$. For a nonnegative integer $n$ and a base $q$, the integer $s_q(n)$ is in fact the minimal number of powers of $q$ needed to represent $n$ as their sum
(which can be proved using that the $q$-ary expansion 
is the lexicographically largest representation of $n$ as a sum of powers of $q$).

Senge and Straus~\cite{SengeStraus1973}
proved the important theorem that for coprime integers $p,q\geq 2$ and arbitrary $c>0$, there are only finitely many integers $n\geq 0$ such that
\begin{equation}\label{eqn_SengeStraus1973}
s_p(n)\le c\quad\text{and}\quad s_q(n)\leq c.
\end{equation}
This statement is, at least heuristically, close to Furstenberg's conjecture~\eqref{eqn_furstenberg}: digital expansions of a number in multiplicatively independent bases usually cannot be simple simultaneously.
Extensions of~\eqref{eqn_SengeStraus1973} were proved by Stewart~\cite{Stewart1980}, Mignotte~\cite{Mignotte1988}, Schlickewei~\cite{Schlickewei1990b,Schlickewei1990}, Peth\H{o}--Tichy~\cite{PethoTichy1993}, and Ziegler~\cite{Ziegler2019}.
See also~\cite{BravoLuca2016,BugeaudCipuMignotte2013,Luca2003} for related results. 

Gelfond~\cite{Gelfond1968} proposed to prove that
\begin{equation}\label{eqn_Gelfond1968}
\#\bigl\{
n\leq x:s_{q_1}(n)\equiv \ell_1\bmod m_1\text{ and }
s_{q_2}(n)\equiv \ell_2\bmod m_2
\bigr\}
=\frac x{m_1m_2}+\LandauO\bigl(x^\delta\bigr)
\end{equation}
for some $\delta<1$,
where $q_1,q_2\geq 2$ are coprime bases, $m_1,m_2$ are integers satisfying $\gcd(m_1,q_1-1)=\gcd(m_2,q_2-1)=1$, and $\ell_1,\ell_2\in\mathbb Z$.
A weak error term $o(1)$ for this problem was proved by B\'esineau~\cite{Besineau1972}, while the full statement was obtained by D.-H.~Kim~\cite{Kim1999}.

Drmota~\cite[Theorem~4]{Drmota2001} proved (among other things) an asymptotic formula for the proportion
\begin{equation}\label{eqn_Drmota2001}
\frac 1x\#\bigl\{n<x:s_{q_1}(n)=k_1,s_{q_2}(n)=k_2\bigr\},
\end{equation}
where $q_1,q_2\geq 2$ are coprime bases,
with an error term $(\log x)^{-1}$.
This may be called a \emph{local limit theorem} for the joint sum-of-digits function $n\mapsto (s_p(n),s_q(n))$.
Note that B\'esineau's result follows as a special case, as the two sum-of-digits functions on $[0,x)$ are mostly found close to their expected values, compare~\eqref{eqn_hoeffding_base2} below.

We also wish to note the recent paper by Drmota, Mauduit, and Rivat~\cite{DrmotaMauduitRivat2019b}, who proved a result on the sum of digits of prime numbers in two different bases.

The starting point for the present paper is the article~\cite{DHLL2017} by Deshouillers, Habsieger, Landreau, and Laishram.
\begin{quote}
``[\ldots] it seems to be unknown whether there are infinitely many integers $n$ for which $s_2(n)=s_3(n)$ or even for which
$\lvert s_2(n)-s_3(n)\rvert$ is significantly small.''~\cite{DHLL2017}
\end{quote}
They prove the following result.
\begin{theorem*}
For sufficiently large $N$, we have
\[\#\bigl\{n\leq N:\lvert s_3(n)-s_2(n)\rvert\leq 0.1457205\log n\bigr\}>N^{0.970359}.\]
\end{theorem*}
Note that the difference $s_3(n)-s_2(n)$ is expected to have a value around
$C\log n$, where
\[C=\frac 1{\log 3}-\frac 1{\log 4}=0.18889\ldots;\]
by the above theorem there exist indeed many integers $n$ such the difference $\lvert s_2(n)-s_3(n)\rvert$ is ``significantly small''.

This result was extended by La Bret\`eche, Stoll, and Tenenbaum~\cite{BretecheStollTenenbaum2019}, who proved in particular that
\begin{equation}\label{eqn_BretecheStollTenenbaum2019}
\bigl\{s_p(n)/s_q(n):n\geq 1\}
\end{equation}
is dense in $\mathbb R^+$ for all multiplicatively independent integer bases $p$, $q\geq 2$.

We also wish to note the papers~\cite{MauduitSarkozy1997} by Mauduit and S\'ark\"ozy,
and by Mauduit, Pomerance, and S\'ark\"ozy~\cite{MauduitPomeranceSarkozy2005}.
In these papers, integers with a fixed sum of digits and corresponding asymptotic formulas are studied, and possible extensions to several bases are addressed.

Let us call a natural number $n$ such that $s_2(n)=s_3(n)$ a \emph{collision} (of $s_2$ and $s_3$).
The question on the infinitude of collisions, mentioned in~\cite{DHLL2017}, is not a new one.
M.~Drmota (private communication to the author) received a hand-written letter from A.~Hildebrand more than twenty years ago, in which the very same problem was presented.

In the present paper, we give a definite answer to this question.
\begin{theorem}\label{thm_collisions}
There exist infinitely many nonnegative integers $n$ such that
\begin{equation}\label{eqn_collision}
s_2(n)=s_3(n).
\end{equation}
More precisely, for all $\delta>0$ we have
\begin{equation}\label{eqn_quantitative_collision}
\#\bigl\{n<N:s_2(n)=s_3(n)\bigr\}\gg N^{\frac{\log3}{\log4}-\delta},
\end{equation}
where the implied constant may depend on $\delta$.
Note that 
$\log3/\log4=0.792\ldots$.
\end{theorem}

The difficulty in proving this theorem lies in the separation of the values of $s_2(n)$ and $s_3(n)$.
The sum-of-digits functions can be thought of as a sum of independent, identically distributed random variables, and they concentrate (according to Hoeffding's inequality, for example) around the values $\frac12\log_2N$ and $\log_3N$ respectively, where $0\leq n<N$.
More precisely, the variances are of order $\log N$, and the tails of these distributions decay as least as fast as $\exp(-C(x-\mu)^2/\sigma^2)$, where $\mu$ is the expected value, and $\sigma^2$ the variance.
Since the gap $(1/\log3-1/\log4)\log N$ comprises $\asymp (\log N)^{1/2}$ standard deviations, we can only expect a number $\ll N^\delta$ of collisions, where $\delta<1$ is some constant.
In the light of this argument, we see that our result cannot be too far from the true number of collisions.

The increasing sequence $\mathfrak s_{2,3}$ of nonnegative integers $n$ such that $s_2(n)=s_3(n)$ is listed as entry \texttt{A037301} in the OEIS~\cite{OEIS}.
The question whether this sequence is infinite had to remain open there.
The first few collisions are as follows:
\[
\begin{array}{rrrrrrrrrrrrrr@{\hskip 0.25mm}l}
n\mbox{ in}&\mbox{binary}&0&1&110&111&1010&1011&1100&1101&10010&10011&10101&100100\\
n\mbox{ in}&\mbox{ternary}
&0&1&20&21&101&102&110&111&200&201&210&1100\\
n\mbox{ in}&\mbox{decimal}
&0&1&6&7&10&11&12&13&18&19&21&36&.
\end{array}
\]
\begin{remarks}
Note the subsequence $(10,11,12,13)$; contiguous subsequences of $\mathbb N$ of length greater than four do not appear in $\mathfrak s_{2,3}$,
since $s_3$ on such a subsequence contains two consecutive up-steps, while $s_2$ decreases or stays constant after one up-step.
We expect that it is possible to extend our proof to arbitrary \emph{patterns} in $\mathfrak s_{2,3}$: for example, we expect that there are infinitely many $n$ such that
\begin{equation}\label{eqn_run}
s_2(n+v)=s_3(n+v)\quad\text{for}\quad v\in\{0,1,2,3\},
\end{equation}
and infinitely many $n$ (the integer $n=13$ is an example) such that
\begin{equation}\label{eqn_run2}
\bigl\{v\in\{0,\ldots,23\}:s_2(n+v)=s_3(n+v)\bigr\}=\{0,5,6,8,23\}.
\end{equation}
More generally, every pattern that appears at all should appear infinitely often in $\mathfrak s_{2,3}$.
To this end, we will have to study certain residue classes modulo $2^k3^\ell$ --- note that for $n\in (2+8\mathbb Z)\cap (1+9\mathbb Z)$, for example,
we have $s_2(n+v)-s_3(n+v)=c$ for some $c$ and all $v\in\{0,1,2,3\}$.
The next step would be to scan these ``candidate residue classes'' for collisions, using our method.
But residue classes of this form are used in our proof anyway, therefore we are optimistic that the main problems have already been overcome.
(Note that also a suitable replacement for Proposition~\ref{prp_aux} below will have to be found. This proposition takes care of the parity restriction $s_3(n+t)-s_3(n)\equiv s_3(t)\bmod 2$.)

We would like to note that our proof of Theorem~\ref{thm_collisions} is not a constructive one. We do not give an algorithm that allows us to find integers $n$ such that $s_2(n)=s_3(n)$.
We leave it as an open problem to find a construction method for such integers $n$.

Also, it is a very interesting open problem to prove that $s_2(p)=s_3(p)$ for infinitely many prime numbers $p$.
We believe that this question is difficult.
This guess is due to the analogy to \emph{missing digit problems}, where
sparse sets $S\subseteq\mathbb N$ (that is, $\#(S\cap[1,N])\ll N^{\delta}$ for some $\delta<1$) of a similar kind are studied;
Maynard~\cite{Maynard2019}, in an important and difficult paper,
could prove that infinitely many primes excluding any given decimal digit exist.
Our set $S=\{n:s_2(n)=s_3(n)\}$ is even less understood than the set of integers in Maynard's result, hence our scepticism.
\end{remarks}
\bigskip\noindent
\textbf{Plan of the paper.}

The main body of the paper concerns the proof of the auxiliary statement, Proposition~\ref{prp_aux} below, which directly leads to the main theorem.
This proof is organized into three main steps, represented by Propositions~\ref{prp_step_1}--\ref{prp_step_3}.
After the statement of these results, in Section~\ref{sec_finishing}, we prove Proposition~\ref{prp_aux} and thus Theorem~\ref{thm_collisions} from these three propositions.
The three sections thereafter, Sections~\ref{sec_constant},~\ref{sec_small}, and~\ref{sec_critical}, are dedicated to the proofs of the three main steps.
At the end of the paper, we present (mostly difficult) research questions.
\begin{notation}
The symbol $\log$ denotes the natural logarithm,
and $\log_a=\frac 1{\log a}\log$ is the logarithm in base $a>1$.
We use Landau notation, employing the symbols $\LandauO$, $\ll$, and $o$.
The symbol $f(n)\asymp g(n)$ abbreviates the statement
$\bigl(f(n)\ll g(n)$ \textsf{and} $g(n)\ll f(n)\bigr)$,
while $f(n)\sim g(n)$ means that $f(n)/g(n)$ converges to $1$ as $n\rightarrow\infty$.
We also use the exponential $\e(x)=\exp(2\pi i x)$.
For $M\geq 0$, the statement ``$a$ is $M$-close to $b$'' means $\lvert a-b\rvert\leq M$.
\end{notation}
\section{Proofs}
Our main theorem follows from the following proposition.
\begin{proposition}\label{prp_aux}
For all $\delta>0$ the number of $n<N$ such that
\begin{equation}\label{eqn_two_possibilities}
s_2(n)-s_3(n)\in\{0,1\}
\end{equation}
is bounded below by $CN^{\frac{\log3}{\log4}-\delta}$ (where the constant $C$ may depend on $\delta$). 
\end{proposition}
We call an integer $n$ such that~\eqref{eqn_two_possibilities} is satisfied an \emph{almost-collision}.

Let $N\geq 4$ be an integer. We are going to find many collisions in the interval $[N,2N)$ for all large enough $N$, which will prove Theorem~\ref{thm_collisions}.
Let $\varepsilon>0$ be arbitrary throughout this proof.
This variable is used as exponent of $\log\log N$,
and its value, as long as it is strictly positive, is irrelevant for our proof.
For given $N$, we define $\lambda$, $\eta$, $f$, $m$, and $J$ as follows. Set
\begin{equation}\label{eqn_choice}
\begin{array}{r@{\hskip 1mm}lr@{\hskip 1mm}lr@{\hskip 1mm}lr@{\hskip 1mm}lr@{\hskip 1mm}l}
\lambda_0&\coloneqq\log N,&
\eta_0&\coloneqq\lambda_0^{3/4},&
f_0&\coloneqq(\log\lambda_0)^{1/2+\varepsilon},&
m_0&\coloneqq\lambda_0^{1/2}/f_0,&
J_0&\coloneqq f_0^2,\\[1mm]
\lambda&\coloneqq\lfloor\lambda_0\rfloor,&
\eta&\coloneqq4\lfloor\eta_0/4\rfloor,&
f&\coloneqq\lfloor f_0\rfloor,&
m&\coloneqq\lfloor m_0\rfloor,&
J&\coloneqq\lfloor J_0\rfloor.
\end{array}
\end{equation}
We wish to give a rough and very imprecise idea of the meaning of this choice of variables.
The length of a binary or ternary expansion of $n\in[N,2N)$ is of size $\asymp\lambda$, and the standard deviation of a (binary or ternary) sum-of-digits function on $[N,2N)$ is of order $\asymp (\log N)^{1/2}$.
The variable $m$ is smaller than the standard deviation by a factor $f$ (the \emph{fineness}), and taking $J$ steps of length $m$,
we cover sufficiently many standard deviations.
That is, the tail (comprising deviations larger than $Jm$ from the expected value) is bounded by $\lambda^{-D}$ for all $D>0$ due to the presence of $\varepsilon>0$.
Finally, $\eta$ is the ternary length of certain integers $\mathfrak a$ and $\mathfrak b$ that we choose freely.
It is large enough to allow for differences of ternary sum-of-digits functions larger than the standard deviation $\asymp \lambda^{1/2}$ by any logarithmic factor $(\log\lambda)^\rho$
(compare to~\eqref{eqn_mathfraka_choice}), and small enough so that a concatenation of $2J+1$ ternary expansions of length $\eta$ is still much shorter than $\lambda$.

After this very informal explanation of our choice of parameters,
we give a brief description of the proof.
The search for collisions will consist of three main steps.
\begin{enumerate}
\item ``Preparation'':
find a residue class $A'$ on which $f(n+t)-f(n)$ takes prescribed, constant differences, where $f(n)=s_2(n)-s_3(n)$; \label{item_step_1}
\item ``Rarefaction'': concentrate the values of $f(n)$ into the interval $[-Jm,Jm]$ by
finding a rarefied and truncated arithmetic progression $A''\subset A'$, and considering only integers $n\in A''$; \label{item_step_2}
\item ``Fair share'': select only those $n\in A''$ such that $f(n)\in m\mathbb Z$. \label{item_step_3}
\end{enumerate}
Steps~\ref{item_step_2} and~\ref{item_step_3} are used to find many values of $n$ from a given given residue class such that $f(n)\in Q\coloneqq\{-Jm,(-J+1)m,\ldots,Jm\}$.
The purpose of Step~\ref{item_step_1} is to define \emph{in advance} a larger residue class $A'=L+2^\nu3^\beta\mathbb Z$ and a set $\mathbf d=\{d_{-J},d_{-J+1},\ldots,d_J\}$ of \emph{shifts} such that $f(n+d_j)-f(n)=jm+\xi_j$ for all $n\in A'$, all $j\in \{-J,\ldots,J\}$, and some $\xi_j\in\{0,1\}$.
This procedure yields many $n$ such that $f(n)\in\{0,1\}$, by choosing for each index $n$ such that $f(n)\in Q$ the appropriate shift $d(n)\in\mathbf d$.
A short argument involving differences $s_j(n+1)-s_j(n)$ of sum-of-digits functions on residue classes (where $j\in\{2,3\}$) allows us to get rid of the unpleasant correction term $\xi_j$.

We will prove the following three propositions, corresponding to our three steps.
\begin{proposition}\label{prp_step_1}
Let $\beta=(2J+1)\eta+1$ and choose the integer $\nu\geq 1$ minimal such that $2^{\nu-1}\ge3^\beta$.
Set
\begin{equation}\label{eqn_dj_def}
d_j\coloneqq \bigl(\tL^{(j+1+J)\eta}\tO\bigr)_3=3\frac{3^{(j+1+J)\eta}-1}2.
\end{equation}
There exists $L\in\{0,\ldots,2^\nu3^\beta-1\}$ such that $L\equiv 9\bmod 12$, and $\xi_j\in\{0,1\}$ for $-J\leq j\leq J $ such that
\begin{equation}\label{eqn_duck}
f(n+d_j)-f(n)=jm+\xi_j
\quad\text{for all }j\in\{-J,\ldots,J\}\text{ and all }n\in A'\coloneqq L+2^\nu3^\beta\mathbb N.
\end{equation}
\end{proposition}

\begin{proposition}\label{prp_step_2}
For an integer $\zeta\geq 0$, define
\begin{equation}\label{eqn_App_def}
A''\coloneqq\bigl(L+2^\nu3^{\beta+\zeta}\mathbb N\bigr)\cap[N,2N)
\end{equation}
and
\begin{equation}\label{eqn_k_interval}
I\coloneqq\bigl\{k\in\mathbb N:N\leq L+2^\nu3^{\beta+\zeta}k<2N\bigr\}.
\end{equation}
Here $\nu$, $\beta$, and $L$ are given by Proposition~\ref{prp_step_1}.
For all $D>0$ there exists a constant $C=C(D)$ such that the following statement holds.

\begin{equation}\label{eqn_important_conclusion}
\begin{minipage}{0.8\textwidth}
There exists a sequence $\bigl(\zeta_N\bigr)_{N\geq 4}$ of nonnegative integers 
such that\\ $\zeta_N\sim \log_3(N)\bigl(1-\log 3/\log 4\bigr)$ as $N\rightarrow\infty$, and for all $N$ and all but at most\\[1mm]
$C\lvert I\rvert\lambda^{-D}$ integers $n\in A''$,
the quantity $f(n)$
is $Jm$-close to $0$.
\end{minipage}
\end{equation}

\noindent
Note that $I$ and $A''$ in this statement depend on $\zeta=\zeta_N$, which in turn depends on $N$.
\end{proposition}

\begin{proposition}\label{prp_step_3}
Using the set $A''$ from~\eqref{eqn_App_def}, we set
\begin{equation}\label{eqn_P_def}
P\coloneqq\#\bigl\{n\in A'':f(n)\in m\mathbb Z\bigr\}.
\end{equation}
As $N\rightarrow\infty$, we have
\begin{equation}\label{eqn_good_class}
P=\frac{\lvert I\rvert}m\bigl(1+o(1)\bigr).
\end{equation}
That is, the residue class $m\mathbb Z$ receives the expected ratio $\lambda^{-1/2}(\log\lambda)^{1/2+\varepsilon}$ of the values of $f(n)=s_2(n)-s_3(n)$ along the finite arithmetic progression $A''$ defined in~\eqref{eqn_App_def}.
\end{proposition}
\subsection{Deriving Theorem~\ref{thm_collisions} from Propositions~\ref{prp_step_1}--\ref{prp_step_3}}\label{sec_finishing}
The expected number $P$ of integers $n\in A''$ such that $f(n)\in m\mathbb Z$ is given by Proposition~\ref{prp_step_3}.
At the same time, Proposition~\ref{prp_step_2} states that for all $D>0$,
$f(n)$ lies in the interval $[-Jm,Jm]$ for $\lvert I\rvert(1-\LandauO(\lambda^{-D}))$ many integers $n\in A''$ (where the implied constant depends on $D$).
Note that for $D>1/2$ this error term is of smaller magnitude than $P$.
Consequently, any choice $D>1/2$ will yield many integers $n\in A''$ such that
$s_2(n)-s_3(n)=jm$ for some $j\in\{-J,\ldots,J\}$.
By~\eqref{eqn_duck} the integer $n'=n+d_{-j}$ satisfies
$s_2(n')-s_3(n')\in\{0,1\}$.
Noting that $\zeta\asymp \log N$ and $J\eta\ll (\log N)^{3/4}(\log\log N)^{1+2\varepsilon}$, we see that the shifts $d_j$ are asymptotically smaller than the common difference $2^\nu 3^{\beta+\zeta}$ of $A''$.
Varying $N$, we get an almost-collision (as in Proposition~\ref{prp_aux}) in each large enough interval $[N,2N)$ and thus the qualitative statement in Theorem~\ref{thm_collisions}.

Considering the asymptotic sizes of $\nu$, $\beta$, and $\zeta$, it is easy to see that the interval $I$ defined in~\eqref{eqn_k_interval} is in fact of size
$\gg N^{\log 3/\log 4-\delta}$ for all $\delta>0$.
Most $k\in I$ yield a value $\tilde f(k)=f(L+2^\nu3^{\beta+\zeta}k)\in [-Jm,Jm]$ by~\eqref{eqn_important_conclusion},
and the expected proportion $\sim m^{-1}\gg (\log N)^{-1/2}$ of them satisfy $\tilde f(k)\in m\mathbb Z$, see~\eqref{eqn_good_class}.
These $k$ yield pairwise different values $k+d_{j(k)}$ as before.
Here the integer $j=j(k)$ is chosen suitably from $\{-J,\ldots,J\}$ in order to 
force an almost-collision.

Let $\delta>0$ be given and set $A=\log3/\log4-\delta$.
If the number of $n<N$ such that $s_2(n)-s_3(n)=0$ and $n\equiv 9\bmod 12$
is $\gg N^A$, there is nothing to be done.
Otherwise, we note that $n\equiv 9\bmod 12$ is equivalent to $\bigl(n\equiv 0\bmod 3$ \textsf{and} $n\equiv 1\bmod 4\bigr)$, therefore $s_3(n+1)=s_3(n)+1$ and $s_2(n+1)=s_2(n)$.
The existence of a number $\gg N^A$ of solutions of $s_2(n)-s_3(n)=1$ on $(9+12\mathbb Z)\cap[N,2N)$ therefore implies a number $\gg N^A$ of collisions on $(10+12\mathbb Z)\cap[N,2N)$.
This establishes~\eqref{eqn_quantitative_collision} and completes the proof.
\qed
\begin{remark}
In the last step towards finding almost-collisions --- namely, choosing $j\in\{-J,\ldots,J\}$ suitably --- the ``element of non-constructiveness'' in our argument is clearly visible.
Currently we do not have any control over the choice of $j$.
\end{remark}
In order to prove Theorem~\ref{thm_collisions}, it is sufficient to establish Propositions~\ref{prp_step_1}--\ref{prp_step_3}.
\subsection{Constant differences of sum-of-digits functions --- proof of Proposition~\ref{prp_step_1}}\label{sec_constant}
We will use \emph{blocks} in ternary, whose lengths are given by the integer $\eta$.
Let us choose nonnegative integers $d_{-J},d_{-J+1},\ldots,d_J$ by concatenating such blocks of ternary digits.
Set
\[\mathfrak b\coloneqq \bigl(\tL^\eta\bigr)_3=\frac{3^\eta-1}2,\]
where $\tL^\eta$ denotes $\eta$-fold repetition of the digit $\tL$.
Define $d_j$, for $-J\leq j\leq J$, by $(j+1+J)$-fold concatenation of $\tL^\eta$, with $\tO$ appended at the right, as in~\eqref{eqn_dj_def}.
The emphasis on ``blocks of length $\eta$'' will become clear in the construction of the integers $k_j$ further down (see~\eqref{eqn_kj_choice}).
Since the ternary expansion of $d_j$ consists of blocks $\tL\tL\tL\tL$ and ends with $\tO$, we have $d_j\equiv 0\bmod 12$ (note that $4\mid(\tL\tL\tL\tL)_3=40$).
Choose the integer $\nu\geq 1$ minimal so that 
\begin{equation}\label{eqn_nu_choice}
2^{\nu-1}\geq 3^{(2J+1)\eta+1}.
\end{equation}
In particular,
\begin{equation}\label{eqn_dj_bound}
d_j<2^{\nu-1}.
\end{equation}
The next important step consists in choosing a certain integer
$a\in\{1,\ldots,2^{\nu-1}-1\}$;
its meaning will become clear in a moment.
The size restrictions imply $d_j+a<2^\nu$ for all $j\in\{-J,\ldots,J\}$.
This means in particular that no carry from the $(\nu-1)$th to the $\nu$th digit occurs in the addition $d_j+a$,
which implies the simple but important identity
\begin{align}\label{eqn_diff_expression}
s_2\bigl(2^\nu n+a+d_j\bigr)-s_2\bigl(2^\nu n+a\bigr)
&=
s_2^{(\nu)}(a+d_j)-s_2^{(\nu)}(a)
\end{align}
for all $n\geq 0$.
The function defined by $s_2^{(\nu)}(n)=s_2(n\bmod 2^\nu)$ is the \emph{truncated binary sum-of-digits function}.
Note that the right hand side of~\eqref{eqn_diff_expression} is independent of $n$; we want to use Chebychev's inequality for choosing a value $a$ such that these values are small for all $j\in\{-J,\ldots,J\}$.
In order to obtain an estimate for the variance, needed for Chebychev's inequality, we adapt parts from~\cite{S2020b}.
For integers $t,L\geq 0$ and $j$, we define a probability mass function $\varphi(\bl,t,L)$ by
\begin{equation}\label{eqn_varphi_def}
\varphi(j,t,L)\coloneqq\frac 1{2^L}\#\bigl\{0\leq n<2^L:s^{(L)}_2(n+t)-s^{(L)}_2(n)=j\bigr\},
\end{equation}
and the characteristic function
\begin{equation}\label{eqn_CF_def}
\omega_t(\vartheta,L)\coloneqq\sum_{j\in\mathbb Z}\varphi(j,t,L)\e(j\vartheta)
=\frac 1{2^L}\sum_{0\leq n<2^L}
\e\bigl(\vartheta s_2^{(L)}(n+t)
-\vartheta s_2^{(L)}(n)
\bigr),
\end{equation}
where $\e(x)=\exp(2\pi ix)$.
Noting that
\begin{equation}\label{eqn_trunc_rec}
s_2^{(L+1)}(2n)=s_2^{(L)}(n)\quad\text{and}\quad
s_2^{(L+1)}(2n+1)=s_2^{(L)}(n)+1,
\end{equation}
the proof of the following statement is not difficult and left to the reader.
\begin{lemma}\label{lem_omega}
For all $t,L\geq 0$ and $j\in\mathbb Z$ we have
\begin{equation}\label{eqn_varphi_rec}
\begin{aligned}
\varphi(j,1,L)&=
\begin{cases}
2^{j-2},&-L+2\leq j\leq 1;\\
2^{-L},&j=-L;\\
0,&\text{otherwise,}
\end{cases}\\
\varphi(j,2t,L+1)&=
\varphi(j,t,L),\\
\varphi(j,2t+1,L+1)&=
\frac12\varphi(j-1,t,L)+\frac 12\varphi(j+1,t+1,L).
\end{aligned}
\end{equation}
The characteristic function satisfies
\begin{equation}\label{eqn_omega_rec}
\begin{aligned}
\bigl\lvert\omega_t(\vartheta,L)\bigr\rvert&\leq 1,\\
\omega_{2t}(\vartheta,L+1)&=\omega_t(\vartheta,L),\\
\omega_{2t+1}(\vartheta,L+1)&=\frac{\e(\vartheta)}2\omega_t(\vartheta,L)+\frac{\e(-\vartheta)}2\omega_{t+1}(\vartheta,L)\quad\text{for } t\geq 1.
\end{aligned}
\end{equation}
\end{lemma}

The recurrence~\eqref{eqn_omega_rec} leads to a recurrence for the moments
\begin{equation}\label{eqn_mk_def}
m_k(t,L)\coloneqq \sum_{j\in\mathbb Z}\varphi(j,t,L)j^k
\end{equation}
of $\varphi(\bl,t,L)$.
Using the identity
\begin{equation}\label{eqn_omega_moments}
\omega_t(\vartheta,L)=\sum_{j\in\mathbb Z}\delta(j,t)\e(jx)
=\sum_{k\geq 0}\frac{m_k(t,L)}{k!}\bigl(2\pi i \vartheta\bigr)^k
\end{equation}
(all involved series are absolutely convergent),
we obtain
\begin{equation}\label{eqn_moment_coeff}
m_k(t,L)=\frac{k!}{(2\pi i)^k}\bigl[\vartheta^k\bigr]\omega_t(\vartheta,L),
\end{equation}
from which we can iteratively obtain recurrences for the moments $m_k(t,L)$.

From~\eqref{eqn_varphi_rec} we clearly see that $m_0(t,L)=\sum_{j\in\mathbb Z}\varphi(j,t,L)=1$, $m_1(t,L)=\sum_{j\in \mathbb Z}j\varphi(j,t,L)=0$, $\varphi(j,2t,L+1)=\varphi(j,t,L)$, and $m_2(1,L)=2-2^{-L+1}$.
Moreover,~\eqref{eqn_omega_rec}, ~\eqref{eqn_omega_moments}, and~\eqref{eqn_moment_coeff} imply
\begin{align*}
m_2(2t+1,L+1)&=
-\frac{1}{4\pi^2}
\bigl[\vartheta^2\bigr]
\Bigl(
\bigl(1+2\pi i\vartheta-2\pi^2\vartheta^2\bigr)
\bigl(1-\bigl(2\pi^2\bigr)m_2(t,L)\bigr)
\\&\hspace{5em}+
\bigl(1-2\pi i\vartheta-2\pi^2\vartheta^2\bigr)
\bigl(1-\bigl(2\pi^2\bigr)m_2(t+1,L)\bigr)
\Bigr)
\\&=m_2(t,L)/2+m_2(t+1,L)/2+1.
\end{align*}
Summarizing, for all $k\geq 0$, $t\geq 1$, and $L\geq 0$, we have
\begin{equation}\label{eqn_moment_rec}
\begin{aligned}
m_0(t,L)&=1,\\
m_1(t,L)&=0,\\
m_2(1,L)&=2-2^{-L+1},\\
m_k(2t,L+1)&=m_k(t,L),\\
m_2(2t+1,L+1)&=\frac{m_2(t,L)+m_2(t+1,L)}2+1.
\end{aligned}
\end{equation}
From the recurrence~\eqref{eqn_omega_rec} for the characteristic function we could easily obtain recurrences for the higher moments too (compare~\cite[(2$\cdot$11)]{S2020b}), but here we only need the first and second moments.
In analogy to Corollary~2.3 in~\cite{SWallner2021} we obtain the following statement.
\begin{quote}
There exists a constant $C$ such that for all integers $B,L\geq 1$, and $t\ge1$ having $B$ blocks of $\tL$s,
\begin{equation*}
m_2(t,L)\leq CB.
\end{equation*}
\end{quote}
However, we only need the following version, which follows directly from~\eqref{eqn_moment_rec}: we have
\begin{equation}\label{eqn_m2_bound}
m_2(t,\nu)\leq 2\nu\quad\text{for all }t,\nu\ge1\text{ such that }t<2^\nu.
\end{equation}
In particular, this holds for $t=d_j$ defined in~\eqref{eqn_dj_def},
and for this estimate we do not need to know what $d_j$ looks like in binary.
We are interested in the differences on the right hand side of~\eqref{eqn_diff_expression}.
By Chebychev's inequality and~\eqref{eqn_m2_bound}, the number of integers $a\in\{0,\ldots,2^\nu-1\}$ such that
\begin{equation}\label{eqn_good}
\Bigl\lvert s_2^{(\nu)}(a+d_j)-s_2^{(\nu)}(a)\Bigr\rvert\leq R_2(2\nu)^{1/2}
\end{equation}
is bounded below by
\[
2^\nu\bigl(1-1/R_2^2\bigr).
\]
Intersecting $2J+1$ sets, we obtain the set of $a<2^\nu$ that satisfy~\eqref{eqn_good} for all $j\in\{-J,\ldots,J\}$,
having cardinality $\ge2^\nu(1-(2J+1)/R_2^2)$.
We choose $R_2=\lambda/(2\nu)$, which is $\asymp\lambda^{1/8}/(\log\lambda)^{1/2+\varepsilon}$ as $N\rightarrow\infty$.
It follows that the set of $a\in\{0,\ldots,2^\nu-1\}$ satisfying
\begin{equation}\label{eqn_good2}
\Bigl\lvert s_2^{(\nu)}(a+d_j)-s_2^{(\nu)}(a)\Bigr\rvert\leq\lambda^{1/2}
\quad\text{for all } j\in\{-J,\ldots,J\}
\end{equation}
has at least
\[
2^\nu\bigl(1-\LandauO\bigl((\log \lambda)^{2+4\varepsilon}\lambda^{-1/4}\bigr)\bigr)
\]
elements,
by the definitions~\eqref{eqn_choice}.
Since powers win against logarithms for large $N$,
we obtain some integer $a$ with the properties that
\begin{equation}\label{eqn_good3}
\begin{array}{ll}
a\equiv 1\bmod 4,\\[2mm]
0\leq a<2^{\nu-1},&\text{and}\\[2mm]
\lvert \delta_j\rvert\leq\lambda^{1/2}
&\text{for all }j\in\{-J,\ldots,J\},
\end{array}
\end{equation}
where
\begin{equation}\label{eqn_delta_def}
\delta_j\coloneqq s_2^{(\nu)}(a+d_j)-s_2^{(\nu)}(a).
\end{equation}

Note that the first two restrictions in~\eqref{eqn_good3} will pose no problem since asymptotically almost all $a<2^\nu$ (as $N\rightarrow\infty$) satisfy the third.

By~\eqref{eqn_diff_expression} we have therefore found an arithmetic progression
\begin{equation}\label{eqn_AP_def}
A=a+2^\nu\mathbb N
\end{equation}
such that each of the sequences
\[\sigma_j=\bigl(s_2(m+d_j)-s_2(m)\bigr)_{m\in A},\]
for $-J\leq j\leq J$, is constant, and attains a value $\delta_j$ bounded by $\lambda^{1/2}$ in absolute value.

In the next step, the ternary sum of digits will come into play,
and we rarefy the progression $A$ by a factor $3^\beta$,
where
\begin{equation}\label{eqn_beta_def}
\beta=(2J+1)\eta+1.
\end{equation}
Note that $\eta\asymp \lambda^{3/4}$ has been used in the definition~\eqref{eqn_dj_def} of the values $d_j$ before.
The selection of this subsequence has to be carried out with care,
so that certain differences $f(n+d_j)-f(n)$, where
\begin{equation}\label{eqn_f_def}
f(n)=s_2(n)-s_3(n),
\end{equation}
are attained on this rarefied progression for $-J\leq j\leq J$.
Sure enough, in order to obtain these differences we will have to ``repair'' the deviation $\delta_j$ from $0$ caused by the differences of binary sums of digits.
We are going to select a residue class $B=K+3^\beta\mathbb N$, where $K<3^\beta$,
on which certain differences
\begin{equation}\label{eqn_base3_differences}
s_3(n+d_j)-s_3(n)
\end{equation}
occur for $n\in B$.
This process will be executed step by step,
thinning out the current residue class by a factor $3^\eta$ for each $j\in\{-J,\ldots,J\}$.
We have found a certain arithmetic progression $A$ in~\eqref{eqn_AP_def}.
A sub-progression $A'$ of $A$ having the desired difference properties in bases $2$ and $3$ --- that is, $s_2(n+d_j)-s_2(n)=\delta_j$ and~\eqref{eqn_diff_expression_ternary} below
---
will be obtained by the intersection
\begin{equation}\label{eqn_intersection}
A\cap B=\left(a+2^\nu\mathbb N\right)\cap\left(K+3^\beta\mathbb N\right)
=L+2^\nu3^\beta\mathbb N,
\end{equation}
where $0\leq L<2^\nu3^\beta$.
We need to find $K$.
This number will in fact be divisible by $3$ (hence the definition of $d_j$ as a multiple of three) ---
together with $a\equiv 1\bmod 4$ this leads to $L\equiv 9\bmod 12$.
The construction is similar to the definition of $d_j$, where we concatenated ternary expansions of length $\eta$, given by $\mathfrak b=\bigl(\tL^{\eta}\bigr)_3$.
We begin with the integer $k_{-J}$.
By our preparation, the quantity $Jm+\delta_{-J}$ (of size $\lambda^{1/2}$ times a logarithmic factor) is considerably smaller than $\eta$ (of size $\lambda^{3/4}$).

The large number of $\tL$s in $\mathfrak b$ can be used to find some $\mathfrak a\in\{0,\ldots,3^{\eta-1}-1\}$ and $\xi\in\{0,1\}$ such that
\begin{equation}\label{eqn_mathfraka_choice}
s_3\bigl(\mathfrak a+\mathfrak b\bigr)-s_3(\mathfrak a)
=
Jm+\delta_{-J}-\xi.
\end{equation}
In fact, such an integer $\mathfrak a$ is found by assembling blocks of length four of ternary digits, where no carry between these blocks occurs,
using the following addition patterns in base $3$:
\begin{equation*}
\begin{array}{r@{\hskip 1mm}r@{\hskip 0em}l}
 &\tO\tL\tL\tZ\\
+&\tL\tL\tL\tL\\
\hline
=&\tZ\tO\tO\tO&.
\end{array}
\quad
\begin{array}{r@{\hskip 1mm}r@{\hskip 0em}l}
 &\tO\tZ\tO\tZ\\
+&\tL\tL\tL\tL\\
\hline
=&\tZ\tO\tZ\tO&,
\end{array}
\quad
\begin{array}{r@{\hskip 1mm}r@{\hskip 0em}l}
 &\tO\tZ\tO\tO\\
+&\tL\tL\tL\tL\\
\hline
=&\tZ\tO\tL\tL&.
\end{array}
\end{equation*}
We see that each block of length four can be used to obtain a variation $\in\{-2,0,2\}$ of the ternary sum of digits;
there are $\eta/4\gg \lambda^{3/4}$ such blocks,
while the needed variation is $\asymp\lambda^{1/2}(\log\lambda)^{1/2+\varepsilon}$
and thus much smaller.
Moreover, by construction~\eqref{eqn_choice}, the integer $\eta$ is divisible by four, so there are no phenomena due to trailing digits.
Using any $\xi\in\{0,1\}$ and $\mathfrak a<3^{\eta-1}$ satisfying~\eqref{eqn_mathfraka_choice},
we set
\begin{equation}\label{eqn_kminusJ_choice}
k_{-J}\coloneqq3\mathfrak a\quad\text{and}\quad\xi_{-J}\coloneqq \xi.
\end{equation}
Trivially, we obtain
\begin{equation}\label{eqn_kminusJ}
s_3\bigl(k_{-J}+d_{-J}\bigr)-s_3\bigl(k_{-J}\bigr)
=
Jm+\delta_{-J}-\xi_{-J}.
\end{equation}
Since $\mathfrak a<3^{\eta-1}$, there does not appear a carry to the $\eta+1$th ternary digit in the addition $k_{-J}+d_{-J}$.
Assume that $k_{j-1}$ has already been defined, for some $-J<j\leq J$.
In analogy to the above, choose $\mathfrak a\in\{0,\ldots,3^{\eta-1}-1\}$ and $\xi\in\{0,1\}$ in such a way that
\begin{equation}\label{eqn_mathfraka_general_choice}
s_3(\mathfrak a+\mathfrak b)-s_3(\mathfrak a)=-m-\delta_{j-1}+\xi_{j-1}+\delta_j-\xi,
\end{equation}
and set
\begin{equation}\label{eqn_kj_choice}
k_j=k_{j-1}+3^{(j+J)\eta+1}\mathfrak a\quad\text{and}\quad\xi_j\coloneqq\xi.
\end{equation}
Note that the target value satisfies $-m-\delta_{j-1}+\xi_{j-1}+\delta_j-\xi\ll\lambda^{1/2}$,
which is again small compared to the number of $\tL$s in $\mathfrak b$.
Since carry propagation between blocks of length $\eta$ is not possible by construction (as in the case $j=-J$), we obtain by concatenating blocks of length $\eta$ and applying a telescoping sum,
\begin{equation}\label{eqn_kj_general}
s_3\bigl(k_j+d_j\bigr)-s_3\bigl(k_j\bigr)
=
-jm+\delta_j-\xi_j
\quad\text{for all }j\in\{-J,\ldots,J\}.
\end{equation}
Finally, set $K=k_J$ and note that
$\beta=(2J+1)\eta+1$ according to~\eqref{eqn_beta_def},
so that $K<3^\beta$.
By construction (note that the ternary digits of $d_j$ from $(j+1+J)\eta+1$ on are zero) we have
\begin{equation}\label{eqn_ternary_part}
s_3(K+d_j)-s_3(K)=-jm+\delta_j-\xi_j
\quad\text{for all }j\in\{-J,\ldots,J\}.
\end{equation}
Similar to~\eqref{eqn_diff_expression}, noting that there is no carry propagation in base three to the $\beta$th digit in the addition $K+d_j$,
we have in fact
\begin{equation}\label{eqn_diff_expression_ternary}
s_3(n+d_j)-s_3(n)=-jm+\delta_j-\xi_j
\end{equation}
for all $n\in K+3^\beta\mathbb N$.
Define $L$ by~\eqref{eqn_intersection}.
By construction, the residue class $L+2^\nu3^\beta\mathbb Z$ is a subset of both $3\mathbb Z$ and $1+4\mathbb Z$, therefore $L\equiv 9\bmod 12$, and we obtain the \emph{difference property}~\eqref{eqn_duck} and thus Proposition~\ref{prp_step_1}.\qed
\subsection{Small values of $f(n)$ --- proof of Proposition~\ref{prp_step_2}}~\label{sec_small}
By our difference property~\eqref{eqn_duck} it is sufficient to prove the existence of (many) elements $n\in A'$ such that
\begin{equation}\label{eqn_f_restriction}
f(n)\in Q,
\quad\text{where}\quad Q=\{jm:-J\leq j\leq J\}.
\end{equation}
After all, for each $n$ satisfying~\eqref{eqn_f_restriction} we can adjust the value of $f$, up to a correction term $\in\{0,1\}$, by any amount $c\in Q$ using a suitably chosen shift $d(n)\in\{d_{-J},d_{-J+1},\ldots,d_J\}$.
Having done so, we arrive at the desired property $f(n+d(n))\in\{0,1\}$.
Since for each given $N$ the constructed quantities $d_j$ are nonnegative and smaller than the common difference of $A'$ --- by~\eqref{eqn_dj_def} we have $d_j<2^\nu3^{\beta}$ ---
this will show that there are infinitely many solutions to $s_2(n)-s_3(n)\in\{0,1\}$, and in fact we will give a quantitative lower bound.
Proving that~\eqref{eqn_f_restriction} has many solutions in $A'$ will be the subject of this and the following section,
constituting the second (``rarefaction'') and third (``fair share'') stages of our proof, respectively.

In the present section we are concerned with restricting our residue class $A'$ in order to obtain $f(n)\in[-Jm,Jm]$ for many integers $n$ in the new set $A''$.
The third step will consist in the study of the property $f(n)\in m\mathbb Z$, which will be carried out in Section~\ref{sec_critical}.

Note that for all $M$, the value $s_2(a+nM)$ will be $C\sqrt{\log N}$-close to $\log_4(N)$ for asymptotically almost all $n<N$ as $N\rightarrow\infty$,
while $s_3(a+nM)$ will be $C\sqrt{\log N}$-close to $\log_3(N)$ most of the time.
Therefore a concentration property of $f(n)$ can only be satisfied for a finite segment of any arithmetic progression.
The fact that the values of $f$ can be concentrated around zero by selecting a finite arithmetic subsequence is an essential point.
It is based on the consideration that $3^\tau n$ has the same ternary sum of digits as $n$ for all integers $\tau\geq 0$,
while the binary sum of digits --- usually --- increases considerably under multiplication by $3^\tau$.
This small remark is in fact the main idea that started the research on the present paper.

Recall the definition~\eqref{eqn_App_def} of $A''$,
for a natural number $\zeta$ that will be chosen in due course.
Suitable choice of $\zeta$ will cause most values of $f$ along $A''$ to lie in the interval $[-Jm,Jm]$.
At this point we only note that $3^\zeta$ will be much larger than $2^\nu$ and $3^\beta$,
in orders of magnitude, $\nu\asymp\beta\asymp \lambda^{3/4}(\log\lambda)^{1+2\varepsilon}$, while $\zeta\asymp \lambda$.
Trivially,~\eqref{eqn_duck} is satisfied on the subsequence $A''$ too.
We are therefore interested in the expression
\begin{equation}\label{eqn_critical}
f\bigl(L+2^\nu3^{\beta+\zeta}k\bigr)=
s_2\bigl(L+2^\nu3^{\beta+\zeta}k\bigr)
-s_3\bigl(L+2^\nu3^{\beta+\zeta}k\bigr),
\end{equation}
where $k$ varies in the interval $I$ defined in~\eqref{eqn_k_interval}.
We can decompose~\eqref{eqn_critical} in the form
\begin{equation}\label{eqn_critical_decomposition}
f\bigl(L+2^\nu3^{\beta+\zeta}k\bigr)=
s_2\bigl(b_2+3^{\beta+\zeta}k\bigr)-s_3\bigl(b_3+2^\nu k\bigr)
+s_2(r_2)-s_3(r_3),
\end{equation}
where
\begin{equation*}
\begin{array}{r@{\hskip 2mm}c@{\hskip 2mm}l@{\hskip 2em}l@{\hskip 2em}r@{\hskip 2mm}c@{\hskip 2mm}l}
b_2&=&\bigl\lfloor 2^{-\nu}L\bigr\rfloor&\text{and}&
b_3&=&\bigl\lfloor 3^{-\beta-\zeta}L\bigr\rfloor,\\[1mm]
r_2&=&L\bmod 2^\nu&\text{and}&r_3&=&L\bmod 3^{\beta+\zeta}.
\end{array}
\end{equation*}
Let us choose
\begin{equation}\label{eqn_zeta_choice}
\zeta_0\coloneqq\log_3(N)\left(1-\frac{\log3}{\log4}\right)+
s_3(L)-s_2(r_2)+\frac\nu2-\beta,
\quad\text{and}\quad
\zeta\coloneqq\lfloor\zeta_0\rfloor.
\end{equation}
We have $r_2<2^\nu$, and $L<2^\nu3^\beta$;
moreover, it follows from the definitions that $\nu=o(\log N)$ and $\beta=o(\log N)$.
Therefore $\zeta\sim C\log_3N$, where the constant equals
\begin{equation}\label{eqn_C_value}
C=1-\frac{\log3}{2\log2}=0.207\ldots.
\end{equation}
In particular, $3^\zeta\geq 2^\nu$ for all large $N$.
Since $L<2^\nu3^\beta$, we have in fact
\[b_3=0\quad\mbox{and}\quad r_3=L.\]
That is, $r_2$ and $r_3$ do not depend on the particular choice of $\zeta\geq \nu\log_32$.
In~\eqref{eqn_zeta_choice} this freedom is used in order to define the rarefaction parameter $\zeta$ suitably.
This in turn determines the arithmetic progression $A''$ defined in~\eqref{eqn_App_def}.
Note that we have already replaced $r_3$ by $L$ in the definition of $\zeta_0$ in order to avoid a circular definition.
This procedure, as we will see, very accurately defines an interval around zero in which $f(n)$, for $n\in A''$, can be found most of the time.
That is,~\eqref{eqn_critical} is close to zero for most $k\in I$.

We study the values
\begin{equation}\label{eqn_f2}
f_2(k)=s_2\bigl(b_2+3^{\beta+\zeta}k\bigr)
\quad\text{and}\quad
f_3(k)=s_3\bigl(2^\nu k\bigr)
\end{equation}
separately, as $k$ varies in $I$.

Sure enough, the study of~\eqref{eqn_f2} will be infeasible in general using current techniques.
This is the case because we encounter problems arising from powers of $2$ and $3$,
as considered in the introduction.
In our application however, the interval $I$ is of the form
\begin{equation}\label{eqn_interval_form}
I=[M,2M+\LandauO(1)]
\end{equation}
for some $M$ considerably larger than $2^\nu$ and $3^{\beta+\zeta}$,
which enables us to prove a nontrivial statement on the distributions of $f_2(k)$ and $f_3(k)$.

In the following, we use the abbreviation $\alpha=\beta+\zeta$.
Let us partition the binary expansion of $b_2+3^\alpha k$ into two parts,
using the integer $\kappa_2=\min\{m:2^m\ge3^\alpha\}$.
For all integers $k\geq 0$, we have
\begin{align}\label{eqn_s2_split}
s_2\bigl(b_2+3^\alpha k\bigr)
=s_2\left(\left\lfloor k\frac{3^\alpha}{2^{\kappa_2}}+\sigma\right\rfloor\right)+s_2\bigl(\bigl(b_2+3^\alpha k\bigr)\bmod 2^{\kappa_2}\bigr),
\end{align}
where $\sigma=b_22^{-\kappa_2}<1$, which follows from
$b_2\leq L2^{-\nu}<3^\beta<2^{\kappa_2}$.

The values of $\lfloor k3^\alpha/2^{\kappa_2}+\sigma\rfloor$ start at
$\tilde M+\LandauO(1)$, where
$\tilde M=\rho M$ and $\rho=3^\alpha/2^{\kappa_2}\in (1/2,1)$,
increase step by step as $k$ runs through $I$,
and remain on the same integer for at most two consecutive values of $k$.
Consequently, the distribution of the first summand for $k\in I$ originates from the distribution of $s_2(k')$ for $k'\in I'$, where
\[I'=[\tilde M-1,2\tilde M+1],\]
and each number of occurrences is multiplied by a value $\in\{0,1,2\}$.
Therefore, using the binomial distribution,
the first summand in~\eqref{eqn_s2_split} can be found within a short interval containing $\tfrac 12\log_2 M$ most of the time.
More precisely, we apply Hoeffding's inequality.
Construing the binary sum-of-digits function on $[0,2^K)$ as a sum of independent random variables with mean $1/2$,
we obtain for all integers $T\geq 0$ and real $t\geq 0$
\begin{equation}\label{eqn_hoeffding_base2}
\frac1{2^T}
\bigl\{
0\leq n<2^T:
\bigl\lvert s_2(n)-T/2\bigl\rvert\geq t
\bigr\}
\le2\exp\left(-2t^2/T\right).
\end{equation}
We apply this for $t=Jm/5$ and $T$ minimal such that
$2^T\geq 2\tilde M+1$. Note that
\[T\sim\log_2\left(\frac{N}{2^\nu3^{\beta+\zeta}}\right)\asymp \lambda.\]
Note that we used the definition of $\zeta$ for the latter asymptotics.
From~\eqref{eqn_hoeffding_base2} we obtain
\begin{equation}\label{eqn_base2_hoeffding_consequence_part1}
\begin{aligned}
\left\{
k\in I:
\bigl\lvert s_2\left(\lfloor k3^\alpha/2^{\kappa_2}+\sigma\rfloor\right)
-T/2\bigr\rvert\geq t
\right\}
&\leq
2
\left\{
k'\in I':
\bigl\lvert s_2(k')
-T/2\bigr\rvert\geq t
\right\}
\\&\leq
2
\left\{
0\leq k'<2^T:
\bigl\lvert s_2(k')
-T/2\bigr\rvert\geq t
\right\}
\\&\ll
\exp\bigl(-2\lambda(\log\lambda)^{1+2\varepsilon}/(25T)\bigr)
\\&\ll\exp\bigl(-C(\log\lambda)^{1+2\varepsilon}\bigr)
\\&\ll
\lambda^{-D}
\end{aligned}
\end{equation}
for all $D>0$ and some $C$, as $N\rightarrow\infty$.
Meanwhile, the second summand in~\eqref{eqn_s2_split} also follows a binomial distribution, with mean $\kappa_2/2$ and a corresponding concentration property.
For this, it is important to note that the sum over $k$ is longer than $2^{\kappa_2}$ (for large $N$): this is due to the observation, given in~\eqref{eqn_C_value}, that $C<1/2$. Therefore, multiples of the odd integer $3^\alpha$ traverse each residue class modulo $2^{\kappa_2}$ in a uniform way.
After forming an intersection, the value of $f_2(k)=s_2(b_2+3^\alpha k)$ is $2Jm/5$-close to the value
\[
E_2=
\frac 12\log_2\left(\frac N{2^\nu3^{\beta+\zeta}}\right)+\frac 12\log_2 3^{\beta+\zeta}
=\frac 12\log_2(N)-\frac\nu2,
\]
for all but $\LandauO\bigl(\lvert I\rvert\lambda^{-D}\bigr)$ integers $k\in I$.
The contribution of $f_3(k)=s_3(2^\nu k)$ can be handled in an analogous fashion.
In this case, the expression $f_3(k)$ is $2Jm/5$-close to the value
\[E_3=
\log_3\left(\frac N{2^\nu3^{\beta+\zeta}}\right)
+\log_3(2^\nu)
=\log_3(N)-\beta-\zeta
\]
for all but $\LandauO\bigl(\lvert I\rvert\lambda^{-D}\bigr)$ integers $k\in I$.
Again, $D>0$ is arbitrary.
Including the term $s_2(r_2)-s_3(r_3)$ from~\eqref{eqn_critical_decomposition} leads to the definition of $\zeta$ in~\eqref{eqn_zeta_choice}.
Joining the preceding statements and~\eqref{eqn_critical_decomposition},
noting that the allowed deviation $Jm$ is not surpassed when adding two times the error $2Jm/5$ and also considering the rounding error coming from the floor function $\zeta=\lfloor\zeta_0\rfloor$, we obtain Proposition~\ref{prp_step_2}.\qed
\subsection{The critical expression modulo $m$ --- proof of Proposition~\ref{prp_step_3}}\label{sec_critical}
The final piece in the puzzle, which we consider before we proceed to the assembly of these pieces,
is the study of the function $f(n)\bmod m=(s_2(n)-s_3(n))\bmod m$ along arithmetic progressions.

We are going to adapt the Mauduit--Rivat method for digital problems~\cite{DrmotaMauduitRivat2009,DrmotaMauduitRivat2019,DrmotaMauduitRivat2019b,DrmotaRivatStoll2008,
MartinMauduitRivat2014,
MartinMauduitRivat2015,
MartinMauduitRivat2019,MartinMauduitRivat2019b,MauduitRivat2009,MauduitRivat2010,MauduitRivat2015,MauduitRivat2018}, also applied in the papers~\cite{DrmotaMorgenbesser2012,MorgenbesserStoll2013,Muellner2017,Mullner2018,AzaiezMkaouarThuswaldner2014,S2020}.
This will be used in order to obtain a statement concerning the number
$P$ defined in~\eqref{eqn_P_def},
\begin{equation*}
\begin{aligned}
P&=\#\bigl\{n\in A'':f(n)\in m\mathbb Z\bigr\}
\\&=\#\left\{k\in I:
s_2\bigl(b_2+3^{\beta+\zeta}k\bigr)-s_3\bigl(2^\nu k\bigr)
\equiv t\bmod m
\right\},
\end{aligned}
\end{equation*}
where $t=s_3(r_3)-s_2(r_2)$ (see~\eqref{eqn_critical_decomposition}).
In order to handle this quantity, it is sufficient to study
\begin{equation}\label{eqn_S0_def}
S_0=S_0(\vartheta)=\sum_{k\in I}\e\bigl(\vartheta
s_2\bigl(b_2+3^{\beta+\zeta}k\bigr)
-\vartheta s_3\bigl(2^\nu k\bigr)\bigr),
\end{equation}
with $\vartheta=\ell/m$, where $\ell\in\{0,\ldots,m-1\}$.
By orthogonality relations,
\begin{align}\label{eqn_orthogonality}
P=\frac{\lvert I\rvert}m+
\frac 1{m}
\sum_{1\leq b<m}
\e\left(-\frac{bt}{m}\right)
S_0\left(\frac bm\right),
\end{align}
and it is sufficient to find an upper bound for $S_0(\vartheta)$.
We apply van der Corput's inequality (for example,~\cite[Lemme~4]{MauduitRivat2010}), where $R\geq 1$ is chosen later:
\begin{multline*}
\lvert S_0\rvert^2\leq \frac{\lvert I\rvert+R-1}R\sum_{-R<r<R}\left(1-\frac{\lvert r\rvert}R\right)
\\\times
\sum_{\substack{k\in I\\k+r\in I}}
\e\Bigl(
\vartheta\bigl(s_2\bigl(b_2+3^{\beta+\zeta}(k+r)\bigr)-s_2\bigl(b_2+3^{\beta+\zeta} k\bigr)\bigr)
-\vartheta \bigl(s_3\bigl(2^\nu(k+r)\bigr)-s_3\bigl(2^\nu k\bigr)\bigr)
\Bigr).
\end{multline*}

Next, we apply a suitable \emph{carry propagation lemma} in order to ``cut off digits'', that is, to replace $s_2$ and $s_3$ by \emph{truncated sum-of-digits functions}:
\begin{align*}
s_2^{(\mu_2)}(n)&=s_2\bigl(n\bmod 2^{\mu_2}\bigr),\\
s_3^{(\mu_3)}(n)&=s_3\bigl(n\bmod 3^{\mu_3}\bigr),
\end{align*}
where $\mu_2,\mu_3\ge0$ are chosen later.
See~\cite[Lemma~4.5]{S2020} for the base-$2$ version used here;
an analogous statement holds for all bases, and we also need the completely analogous base-$3$ variant (the original statement was given in~\cite[Lemme~5]{MauduitRivat2010}, compare also~\cite[Lemme~16]{MauduitRivat2009}).
We discard the condition $n+r\in I$, and join the cases $r$ and $-r$, in order to obtain
\begin{align}\label{eqn_S0_S1}
\lvert S_0\rvert^2
&\leq
\lvert I\rvert^2
\mathcal O\left(
\frac R{\lvert I\rvert}
+\frac {3^{\beta+\zeta}R}{2^{\mu_2}}
+\frac {2^\nu R}{3^{\mu_3}}
\right)
+
\frac {2\,\lvert I\rvert}R\sum_{0\leq r<R}
\bigl\lvert S_1\bigr\rvert,
\end{align}
where
\begin{equation}\label{eqn_S1_def}
\begin{aligned}
\hspace{8em}&\hspace{-8em}
S_1=
\sum_{k\in I}
\e\Bigl(
\vartheta s_2^{(\mu_2)}\bigl(3^{\beta+\zeta}k+b_2+3^{\beta+\zeta}r\bigr)
-
\vartheta s_2^{(\mu_2)}\bigl(3^{\beta+\zeta}k+b_2\bigr)
\\&
-\vartheta s_3^{(\mu_3)}\bigl(2^\nu k+2^\nu r\bigr)
+\vartheta s_3^{(\mu_3)}\bigl(2^\nu k\bigr)\Bigr)
.
\end{aligned}\end{equation}
Note that the lowest $\mu_2$ binary digits of $b_2+3^{\beta+\zeta}k$
and the lowest $\mu_3$ ternary digits of $2^\nu k$ 
are visited \emph{uniformly and independently} --- this is just the Chinese remainder theorem.

We obtain
\begin{equation}\label{eqn_S1_split}
\begin{aligned}
\hspace{4em}&\hspace{-4em}S_1=
\frac{\lvert I\rvert}{2^{\mu_2}3^{\mu_3}}
\sum_{0\leq n_2<2^{\mu_2}}
\e\bigl(\vartheta s_2^{(\mu_2)}\bigl(n_2+3^{\beta+\zeta} r\bigr)-\vartheta s_2^{(\mu_2)}(n_2)\bigr)
\\&
\times
\sum_{0\leq n_3<3^{\mu_3}}
\e\bigl(\vartheta s_3^{(\mu_3)}\bigl(n_3+2^\nu r\bigr)-\vartheta s_3^{(\mu_3)}(n_3)\bigr)
+
\mathcal O\bigl(2^{\mu_2}3^{\mu_3}\bigr).
\end{aligned}
\end{equation}
For this estimate to be relevant, it is important that the number $C$ defined in~\eqref{eqn_C_value} is smaller than $1/2$:
the interval $I$ has length $\asymp N/(2^\nu3^{\beta+\zeta})$,
and we need to run through $2^\nu3^{\beta+\zeta}$ many integers $n\in I$ in order to apply the Chinese remainder theorem.
In contrast, comparing the bases $2$ and $7$, the corresponding constant
\[C_{2,7}\coloneqq 1-\frac{(2-1)\log 7}{(7-1)\log 2}=0.532\ldots\] 
 will already be greater than $1/2$, so new ideas will be needed for bases of ``very different size''.
Meanwhile, adjacent bases $b$ and $b+1$, for example, can certainly be handled by our method;
the sequence of constants $C_{b,b+1}$ decreases to zero as $b\rightarrow\infty$.

It is sufficient to find a nontrivial estimate for the first factor in~\eqref{eqn_S1_split}, concerning the binary expansion.
We are concerned with the correlation (a characteristic function) we had in~\eqref{eqn_CF_def}:
\[
\omega_t(\vartheta,L)=\frac 1{2^L}
\sum_{0\leq n<2^L}
\e\bigl(\vartheta s_2^{(L)}(n+t)
-\vartheta s_2^{(L)}(n)
\bigr).
\]
Reusing the argument leading to~\cite[Lemma~2$\cdot$7]{S2020b},
and Lemma~\ref{lem_omega},
we obtain the following result.
\begin{lemma}\label{lem_muntjak}
Assume that integers $B\geq 0$ and $L,t\geq 1$ are given such that $t$ contains at least $2B+1$ blocks of $\tL$s, and $t<2^L$.
Then for all real $\vartheta$,
\[\bigl\lvert \omega_t(\vartheta,L)\bigr\rvert\leq \left(1-\frac 12\lVert \vartheta\rVert^2\right)^B.\]
\end{lemma}
Our focus therefore lies on the number $B$ of blocks of $\tL$s in the binary expansion of $3^{\beta+\zeta}r$.
The only thing we need to know about powers of three in this context is the fact that they are odd integers --- we exploit in an essential way the summation over $r$ instead.
The parameter $R$ will be a certain power of $N$;
in this way, the expected size of $B$ is $\gg\lambda$.

Note that counting the number of blocks of $\tL$s in binary amounts to counting the number of occurrences of $\tO\tL$ (where the $\tO$ corresponds to the more significant digit), up to an error $\LandauO(1)$.
For simplicity, we only count such occurrences where the digit $\tL$ in the block $\tO\tL$ occurs at an even index.
For example, in the binary expansion $\tL\tO\tL\tL\tO\tL\tL\tO$ the corresponding number is $1$, whereas there exist three blocks of $\tL$s.
This simplification will, on average, give $1/2$ of the actual expected value, which is sufficient for our purposes. 
We are therefore concerned with the number $\#\tL(n)$ of $\tL$s occurring in the base-$4$ expansion of $n$:
the number of integers $0\leq n<4^K$ such that $\#\tL(n)=\ell$ is given by
\[4^K\binom K\ell(1/4)^\ell(3/4)^{K-\ell}.\]
Suppose that we have $R=4^K$.
Note that 
\[r\mapsto r3^{\beta+\zeta} \bmod 4^K\]
is a bijection of the set $\{0,\ldots,4^K-1\}$.
We abbreviate $\alpha=1-\lVert\vartheta\rVert^2/2$, and obtain by Lemma~\ref{lem_muntjak}
\begin{equation*}
\begin{aligned}
S_2\coloneqq
\hspace{4em}&\hspace{-4em}
\sum_{0\leq r<R}
\left\lvert
\frac1{2^{\mu_2}}
\sum_{0\leq n_2<2^{\mu_2}}
\e\bigl(\vartheta s_2^{(\mu_2)}\bigl(n_2+3^{\beta+\zeta} r\bigr)-\vartheta s_2^{(\mu_2)}(n_2)\bigr)
\right\rvert
\\&\leq
\sum_{0\leq \ell\leq K}
\sum_{\substack{0\leq r<4^K\\\#\tL(r)=\ell} }
\alpha^{\frac{\ell-2}2}
\\&=
4^K
\alpha^{-1}
\sum_{0\leq \ell\leq K}
\binom K\ell
(1/4)^\ell(3/4)^{K-\ell}
\alpha^{\ell/2}
\\&=
4^K
\alpha^{-1}
\bigl(\sqrt{\alpha}/4+3/4\bigr)^K.
\end{aligned}
\end{equation*}
Since $\sqrt{1+x}\leq 1+x/2$ for $x\geq -1$, we have
\begin{equation}\label{eqn_alpha_estimate}
\sqrt\alpha=\bigl(1-\lVert\vartheta\rVert^2/2\bigr)^{1/2}
\leq 1-\frac 14\lVert\vartheta\rVert^2,
\end{equation}
and the inequality
$(1+x)^K=\exp\bigl(K\log(1+x)\bigr)\leq \exp(Kx)$
yields
\begin{equation}\label{eqn_S2_estimate}
S_2\ll4^K
\exp\left(-\frac K{16}\lVert\vartheta\rVert^2\right).
\end{equation}

We translate this back to $S_0$, noting that
$\lVert\vartheta\rVert\geq 1/m\sim\lambda^{-1/2}(\log\lambda)^{1/2+\varepsilon}$:
for some constant $C>0$ (any value $C\in(0,1/16)$ is good enough) we obtain
\begin{align}\label{eqn_residue_classes}
\lvert S_0\rvert^2
&
\ll
\lvert I\rvert^2
\left(
\frac R{\lvert I\rvert}
+\frac {3^{\beta+\zeta}R}{2^{\mu_2}}
+\frac {2^\nu R}{3^{\mu_3}}
+
\exp\bigl(-CK\lambda^{-1}(\log\lambda)^{1+2\varepsilon}\bigr)
\right).
\end{align}
We see that the last term yields a contribution to $S_0$ that is is smaller than the fair share $\lvert I\rvert m^{-1}\sim\lvert I\rvert\lambda^{-1/2}(\log\lambda)^{1/2+\varepsilon}$ as soon as $K\asymp\lambda$, due to the presence of the power $(\log\lambda)^{1+2\varepsilon}$ in the exponent.
For this, we need to choose $R=4^K$ as large as some positive (fixed) power of $N$.
At the same time we have to take care of the other error terms in~\eqref{eqn_residue_classes}.
It is obvious that we can choose $R\asymp N^\iota$, where $\iota$ is small, and $2^{\mu_2}$ resp. $3^{\mu_3}$ larger than $R\,3^{\beta+\zeta}$ resp. $R\,2^\nu$ (by some small power of $N$),
in such a way that $2^{\mu_2}3^{\mu_3}$ is still smaller than $\lvert I\rvert$ (by another power of $N$).
Such a choice is possible by the fact that $\zeta<1/2$,
and we commented on this after~\eqref{eqn_S1_split}.
We therefore obtain~\eqref{eqn_good_class} from~\eqref{eqn_orthogonality} and~\eqref{eqn_residue_classes}, which completes the proof of Proposition~\ref{prp_step_3} and thus the proof of Theorem~\ref{thm_collisions}.\qed
\section{Open problems}
\begin{enumerate}
\item\label{pb_0} Find a construction method for collisions, and for patterns of collisions as in~\eqref{eqn_run},~\eqref{eqn_run2}.
\item\label{pb_1} Prove that there are infinitely many prime numbers $p$ such that
\begin{equation}\label{eqn_pb_2}
s_2(p)=s_3(p).
\end{equation}
\item\label{pb_2}
Prove or disprove the asymptotic formula
\begin{equation}\label{eqn_pb_3}
\#\bigl\{n<N:s_2(n)=s_3(n)\bigr\}\sim cN^\eta
\end{equation}
for some real constants $c$ and $\eta$.
\item \label{pb_3}
Prove an asymptotic formula (in $k$) for the number of solutions of the equation
\begin{equation}\label{eqn_pb_1}
2^{\mu_1}+\cdots+2^{\mu_k}=3^{\nu_1}+\cdots +3^{\nu_k},
\end{equation}
and for the numbers
\[\#\bigl\{n\in \mathbb N:s_2(n)=s_3(n)=k\bigr\}\]
(finiteness in the second case was proved by Senge and Straus~\cite{SengeStraus1973}).
\item
Generalize Theorem~\ref{thm_collisions} and Problems~\ref{pb_0}--\ref{pb_3} to any pair $(q_1,q_2)$ of multiplicatively independent bases, and to arbitrary families $(q_1,\ldots,q_K)$ of pairwise coprime bases $\geq 2$.
It would also be interesting to prove the existence of infinitely many Catalan numbers \emph{exactly divisible} by some power of $a$, where $a\geq 2$ is an arbitrary integer.
This property can be defined by
\begin{equation}\label{eqn_def_exactdivisibility}
a^k\Vert n\Leftrightarrow \left(a^k\mid n\hskip 2mm\mbox{\textsf{and}}\hskip 1mm\gcd(na^{-k},a)=1\right).
\end{equation}
\item Study collisions of integer-valued $k$-\emph{regular sequences}~\cite{AlloucheShallit1992,AlloucheShallit2003} in coprime bases, generalizing the sum-of-digits case.
\end{enumerate}

\subsection*{Acknowledgements.}
The author is grateful to Michael Drmota and Jo\"el Rivat, who introduced him to digital expansions as a research topic, and to Thomas Stoll for proposing related research problems to him.
Moreover, he thanks Jean-Marc Deshouillers for pointing out the article~\cite{DHLL2017}, which was the starting point for the work on the present paper.

\bibliographystyle{siam}
\bibliography{collisions}

\end{document}